# AERODYNAMIC DESIGN AND PERFORMANCE EVALUATION OF PIPE DIFFUSER FOR CENTRIFUGALCOMPRESSOR OF MICRO GAS TURBINE

Sujal Bhavsar[1], Santanu De[1]

[1]Indian Institute of Technology Kanpur, Kanpur, UP, India

## ABSTRACT

*This paper presents a methodology adopted for calculating geometrical parameters of the pipe diffuser of centrifugal compressor of 100 kW Micro Gas Turbine (MGT). The decision was mostly influenced by the parameters of available optimized airfoil type diffuser and considerations on topological constraint of pipe diffuser, while at the same time aided with diffuser maps. The performance evaluation was performed with 3D-RANS based steady CFD simulation using commercially available ANSYS CFX solver. Performance of baseline of pipe diffuser was compared with airfoil diffuser for 100% rotating speed. Detailed analysis and comparison of flow mechanism was done at design mass flow rate (0.806 Kg/s) at 100% rotating speed. The performance evaluation was also computed at different rotating speeds for pipe diffuser and few operating points have been analysed in detail. At last, four different design of pipe diffuser were also taken for analysis to check the effect of few parameters. It was found that at 100% speed, the best pipe diffuser amongst the considered set performed marginally inferior than airfoil type diffuser. With the same pressure ratio, on average, it offers 82.2% total to total isentropic efficiency as opposed to 84.4% in case with airfoil type diffuser. However, little compactness in frontal area can be achieved with pipe diffuser. It was found that two counter-rotating vortices generated due to leading edge geometry of pipe diffuser offers dual effect—beneficiary effect in pseudo- and semi-vaneless space and destabilizing effect in channel space. Under given operating parameters and considered design, latter effect outweighs the former and a large separation at PS dominated by it is found to be inevitable.*

## NOMENCLATURE

| | | |
|---|---|---|
| $2\theta$ | | divergence angle |
| $\alpha$ | | inclination angle |
| $\beta$ | | angular pitch of single pipe |
| $\eta_{t-t}$ | | total to total isentropic efficiency, |
| $\omega$ | | rotational speed in rpm |
| $A$ | | area |
| $C_{po}$ | | total pressure loss coefficient |
| $C_p$ | | pressure recovery coefficient |
| $D$ | | diameter |
| $L$ | | channel length |
| $m$ | | mass flow rate |
| $Ma$ | | mach number |
| $N$ | | number of blades/pipes |
| $P_{r_{t-t}}$ | | total to total pressure ratio |
| $R$ | | radius |
| $R_{3a}$ | | radius at the end of vaneless space |
| $R_{3b}$ | | radius at the intersection of two pipes |
| $R_{tan}$ | | radius of tangential reference circle |
| $B_{th}$ | | throat blockage |
| $Re_D$ | | reynolds number based on diameter |
| $AR$ | | $A_{out}/A_{th}$ |

**SUBSCRIPTS**

| | |
|---|---|
| 1 | impeller inlet |
| 2 | impeller outlet |
| 3 | diffuser inlet |
| 4 | diffuser outlet |
| $t$ | total |
| $th$ | throat |

**ACRONYMS**

| | |
|---|---|
| AR | Area Ratio |
| CEL | CFX Expression Language |
| DOE | Design of Experiment |
| MGT | Micro Gas Turbine |
| MRF | Multi Reference Frame |
| OFAT | One Factor At Time |
| PS | Pressure Side |

SS   Suction Side

# 1 INTRODUCTION

The IIT Kanpur MGT project is aimed at the development of millimeter-scale gas turbine for potential applications like electric generators and propulsion systems. Though the size of the MGT is of the order of few mm, considerably smaller than the conventional engine, the required compressor stage pressure ratio is comparable. It was targeted to build a compressor with a stage $P_{r_{t-t}}$ of 4. The limited amount of diffusion in an impeller calls for a diffuser capable of large pressure recovery with minimum radial distance. At designed mass flow rate, the absolute value of velocity leaving the impeller reaches the value of 385 m/s with $\boldsymbol{Ma}$ approaching close to 1. Hence one additional constraint is to deal with transonic or supersonic inlet flow without causing choking at the throat of the diffuser.

Kenny [1] introduced the concept of pipe diffuser in 1969, who, from experimental evidence, claimed that the pipe diffuser could achieve higher effectiveness than wedge diffuser by 6-9%. Ge Han et. al. [2] performed computational analysis and obtained 1.2% increase in adiabatic efficiency with pipe diffuser over wedge diffuser at design point. Yang Xi et al [3] also obtained a high effectiveness (83.22 %) with an optimized pipe diffuser at design point and recommended the pipe diffuser for transonic/supersonic flow field at impeller outlet. However, at off-design points, wedge diffuser performed better than pipe diffuser [2, 3]. Kenny [4] compared curves of pressure recovery for a pipe diffuser with wedge diffuser for different inlet flow conditions, and found superior performance of the pipe diffuser for a given level of blockage at the throat. The lower throat blockage of the pipe diffuser suggests a smaller geometric throat for the same flow rate and consequently, a more compact diffuser with high effectiveness can be built [4, 5]. On further research on analyzing the reason for better performance of pipe diffuser over wedge diffuser, it was found that the geometrical nature of leading edge of pipe diffuser offers an advantage over conventional wedge diffuser in accommodating non-uniformities and instabilities in inlet flow [4-6]. This reason was further supported by the results of Dean (1973) who claimed that the higher performance in pipe diffuser lies only in entry configuration rather than conical channel geometry compared to flat, straight type [6]. On further investigation on flow mechanism at entry region of pipe diffuser, Ge Han et. al [5] and Zachau et.al [7] obtained two counter-rotating vortices on the PS of the pipe diffuser from computational and experimental analysis, respectively. It seems that the two counter-rotating vortices actually has a beneficiary effect on an overall performance in the case of Kenny [1]. Ge Han e. al [2] and Yang Xi et al. [3] since it is the only unique fluid dynamic phenomena present in the pipe diffuser as opposed to in wedge diffuser. However, it is possible that for particular operating condition, this could lead to higher loss production and substantial separation at PS of the pipe diffuser [7, 8].

Despite such positive commentaries on pipe diffuser, little data about design and performance of pipe diffuser exist in the open literature [2, 6]. Most of the previous studies compared the performance of pipe diffuser only with wedge type diffuser and for specific set of operating conditions. It is therefore the aim of this work to validate the claims made by previous researcher on pipe diffuser and to check its performance against the available optimized airfoil type diffuser for micro turbomachinery application.

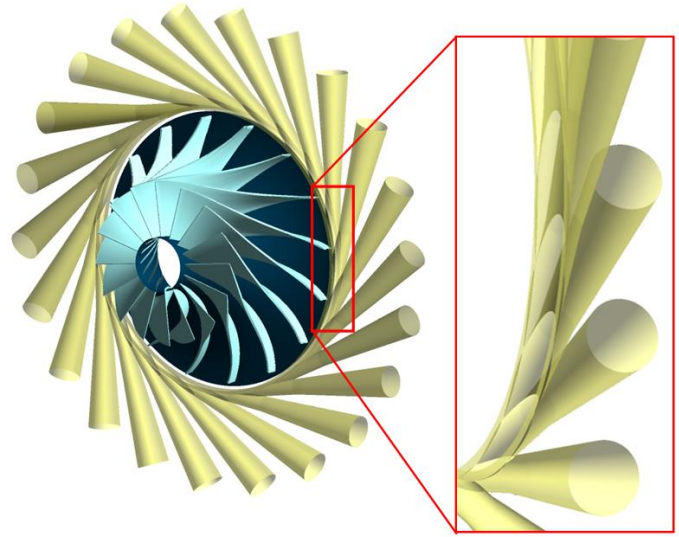

**FIGURE 1:** Pipe diffuser CAD model and view of elliptical leading edge

# 2 DESIGN CONSIDERATION OF PIPE DIFFUSER
## 2.1 Principle difference between geometries of pipe and airfoil diffuser

Vanned type diffuser can be broadly categorized into Channel type and Conical type diffusers [9]. Pipe diffuser falls into the category of latter. Due to intersection of circumferentially arranged pipes, Pipe diffuser features one unique region at the leading region—pseudo-vaneless space, Fig. 2. The other important characteristics of pipe diffuser is the elliptical leading edge, Fig. 1. The other regions such as vaneless, semi-vaneless, and channel space are similar to that of airfoil type diffuser. However, the shape of individual regions is different, Fig. 2.

## 2.2 Geometrical parameters

Comprehensive studies on single straight and conical channel diffuser were available in literature in the form of diffuser maps. It shows the influence of geometrical and fluid-dynamic parameters on the performance of diffuser. We can assume that these maps are also fairly application in the preliminary design stage of the radial diffusers (either Channel type or Conical types) [6]. Francis et al. [10] produced diffuser maps of pressure recovery performance of single conical diffuser as a function of diffuser geometry for fixed set of inlet condition such as $B_t$, $Re_D$, and $Ma_{th}$. Though majority of the decision on geometrical parameters of pipe diffuser were influenced by

2        © 2019 by ASME

airfoil diffuser's parameter and topological constraints, it was also partially supported by the available chart in Francis et al [10].

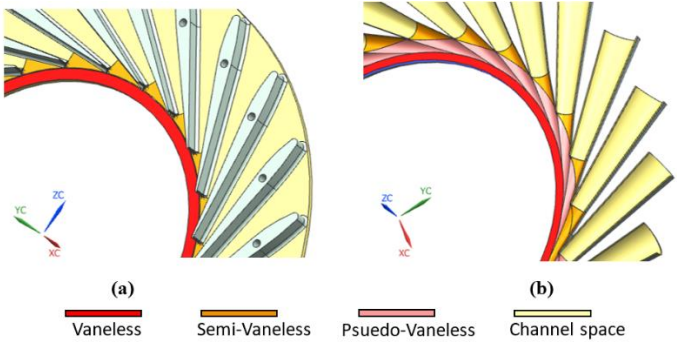

**FIGURE 2:** (a) Airfoil type diffuser (b) Section view of Pipe diffuser

**TABLE 1:** Compressor specifications at design point

| Expected specifications from thermodynamic cycle | |
|---|---|
| $m$ | 4.45 kg/s |
| $P_{r_{t-t}}$ | 4 bar |
| $\eta_{t-t}$ | 84 % |
| Impeller specification | |
| $\omega$ | 70,000 $rpm$ |
| $N$ | 16 |
| $D_2$ | 145 $mm$ |
| Exit blade height | 7.25 $mm$ |
| Optimized airfoil diffuser | |
| $N$ | 23 |
| $AR$ | 4 |
| $L$ | 57 $mm$ |
| $2\theta$ | 10° |
| $\sum A_{th}$ | 817.85 $mm^2$ |
| $R_3/R_2$ | 1.03 |
| $D_4$ | 255 $mm$ |

Ge Han et al [5] concluded from the CFD simulation that constant passage area of the throat region (throat length) leads to local acceleration of flow due to boundary layer thickness growth. Hence, effort is made to make the throat length as minimum as possible (ideally zero). With this consideration, a topological relationship between the constrained parameters was derived and scripted in MATLAB for further use, refer **Appendix A**, which is as:

$$N = f(\alpha, D_{th}, D_{tan})$$

Fig. 3 shows different parameters of pipe diffuser geometry. Table 2 shows group of design parameters for 4 different designs;

**TABLE 2:** Pipe diffuser parameters

| Designs | P1 | P2 | P3 | P4 |
|---|---|---|---|---|
| $N$ | 23 | 22 | 22 | 22 |
| $AR$ | 4 | 4 | 4 | 4 |
| $L\ (mm)$ | 57 | 57 | 45 | 57 |
| $2\theta$ | 6° | 6° | 10° | 6° |
| $\alpha$ | 61.7° | 61.7° | 61.7° | 61.7° |
| $\sum A_{th}$ | 1213.7 | 1154.3 | 1154.3 | 1154.3 |
| $R_{3a}/R_3$ | 1.03 | 1.03 | 1.03 | 1.01 |
| $D_4$ | 255 | 250 | 232 | 250 |
| $D_{tan}$ | 150 | 145 | 145 | 145 |
| $D_{th}$ | 8 | 8 | 8 | 8 |

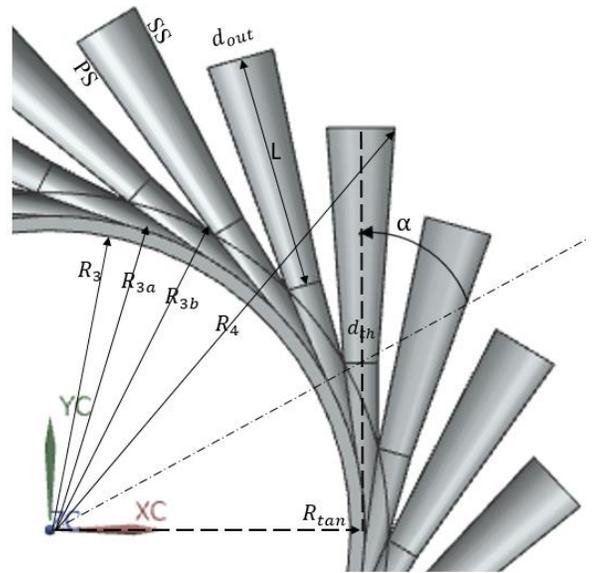

**FIGURE 3:** Geometrical parameters of pipe diffuser

represented as P1, P2, P3, and P4. The effect of $R_{tan}$, $2\theta$, and $R_{3a}/R_3$ was tested. The AR is kept similar to airfoil diffuser– supported by the fact that for isentropic flow, $C_p$ depends only on AR. There is a choice to make either $L$ or $2\theta$ similar. For a similar value of $L$ (57 mm), value of $2\theta$ comes out to be 6°– supported by the diffuser map as the optimal value at high inlet $Ma$ [10]. The alternative consideration was also tested considering in mind that diffuser maps are not absolutely accurate for radial pipe diffuser. In each design iteration, the different combinations of $\alpha$ and $D_{th}$ was tested to obtain a feasible value of $N$. It was found that we don't have much flexibility in increasing $\alpha$ and $D_{th}$; for design P1, increasing $\alpha$ to 66° or increasing $D_{th}$ to 5.5 mm, while keeping other parameters same, would lead to imaginary value of N— indicating failure of design. P1 and P2 were used to check the effect of $R_{tan}$; P2 and P3 were used to check on the choice to make either L or $2\theta$ similar as airfoil diffuser; P2 and P4 were used to check the effect of $R_{3a}/R_3$.

 

## 3 NUMERICAL METHOD AND COMPUTATIONAL PROCEDURE
### 3.1 Simulation setup

The steady aerodynamic simulation performed in this work have been conducted with ANSYS CFX solver. For detailed description, the author referred to the technical guide of ANSYS CFX solver [11]. CFX solver features different class of turbulent models such as RANS, LES, DNS. SST model, which is based on RANS equation, have been used for simulation. SST model was designed to give a highly accurate prediction of the onset and the amount of flow separation under adverse pressure gradients, and is recommended for high accuracy boundary layer simulation [11]. The stage of compressor consists of fast rotating impeller and stationary diffuser. The concept of MRF modeling was adopted to perform steady-state simulation in which each individual cell zone can be assigned different rotating/translating speeds. The RANS equations transform from inertial frame into rotating reference in impeller zone [11]. The interface model named *Stage* (Mixing Plane) was applied between impeller and diffuser zone. *Stage* model performs a circumferential averaging of the fluxes across different frames and/or pitch [11]. Physical Timestep method was used for timescale control. The adequate physical timestep was given via CEL expression.

Inlet condition of impeller is specified as subsonic regime, total pressure (100830 Pa) at designed point, flow direction profile, and 5% turbulent intensity. Outlet condition of diffuser is specified as subsonic flow regime, and mass flow rate (0.806 kg/s at designed point). The exit of the impeller and inlet of the diffuser makeup of the *mixing plane interface*. Periodic boundaries were assigned as shown in Fig. 4. Other boundaries were kept as adiabatic wall. By varying rotational speed and mass flow rate, performance map was generated. Margin on the magnitude of residuals was kept below 1e-4. Surge and choke limit were calculated based on convergence difficulty.

### 3.2 Meshing

The impeller meshing is a structured multi-block hexahedral mesh having 897126 nodes. It was created using ANSYS Turbogrid module which automates the production of high-quality structured mesh needed for blade passage [12]. The diffuser meshing is unstructured tetrahedral mesh created using ANSYS ICEM CFD [13] having 448358 nodes, Fig. 5. In case of pipe diffuser, the surface mesh was first generated using robust *Octree* Method which uses top-down approach, and then regeneration of volume mesh was done using *Delaunay* method which uses bottom-up approach. By following this methodology, the author was able to generated adequate transition in volume meshing. In order to resolve near-wall boundary layer flow in diffuser, the distance away from the wall of the first node was calculated (7e-7 m) such that $y^+$ value becomes 1. To generate adequate inflation layer near wall, 22 prism layers was used. The grid independence test was performed individually on both impeller and diffuser before arriving at the medium quantity of grid.

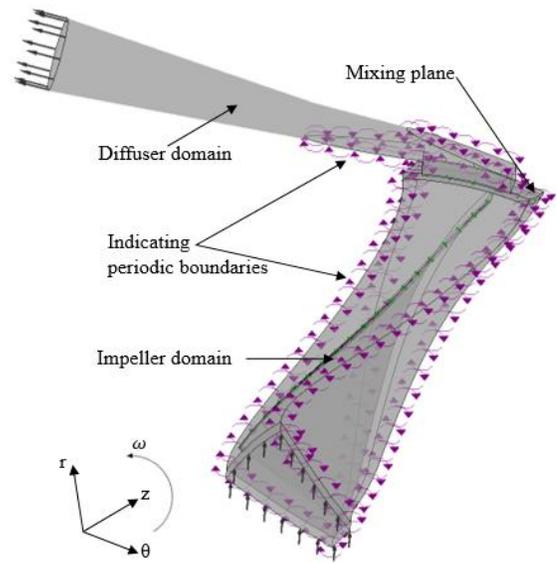

**FIGURE 4:** Boundary conditions in simulation setup

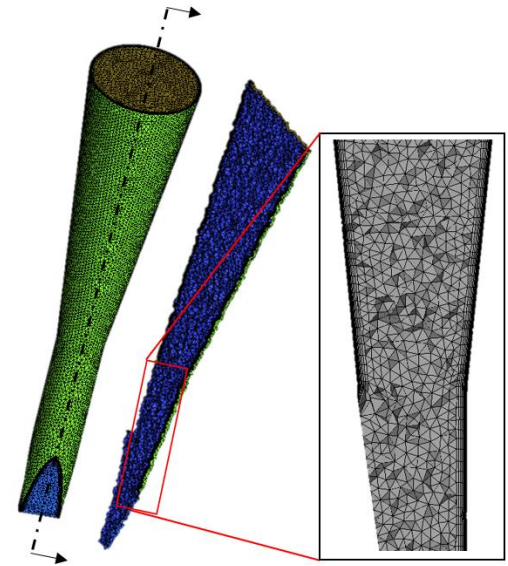

**FIGURE 5:** Full and section view of pipe diffuser mesh

### 3.3 Numerical validation

The developed numerical model was validated with experimental results for airfoil diffuser. Fig. 7 compares the numerically predicted performance map with the experimental results for airfoil diffuser. It gives the variation of pressure ratio with mass flow parameter. Solver over predict the pressure ratio by +2.5% at design point. Surge line is decided by the converged point with maximum 5% highest backpressure. Solver gives good agreement with experimental surge line with maximum under prediction of -0.2%. Choke margin is decided at the point where solver fails. Solver over predict the pressure ratio at low

4© 2019 by ASME

speed and under predict at higher speed. The maximum deviation in the pressure ratio is -5.5% which occur at higher speed near surge line.

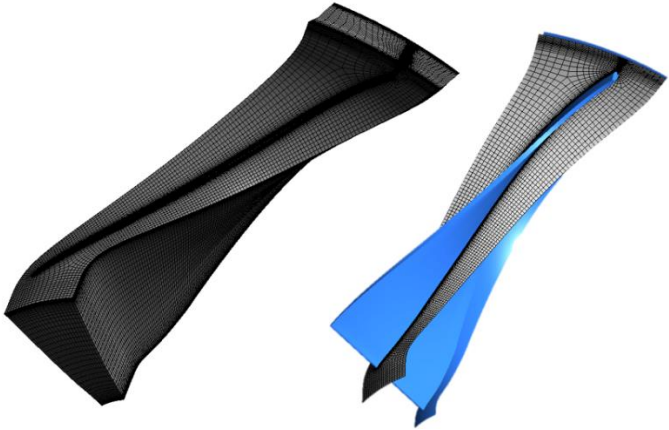

**FIGURE 6:** Full and 50% span view of impeller blade mesh

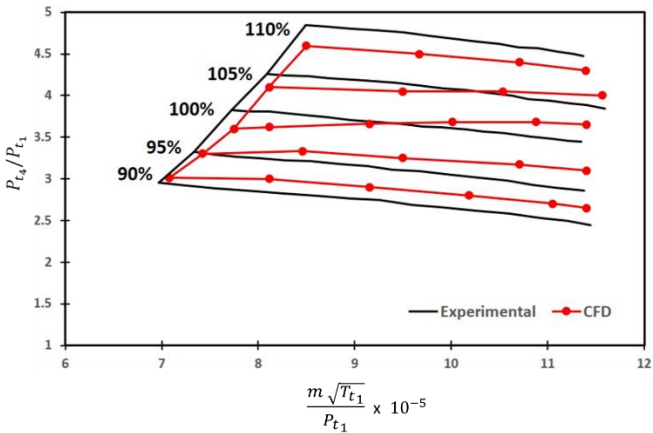

**FIGURE 7:** Comparison of CFD predicted and experimental results for airfoil diffuser

## 4 RESULTS AND DISCUSSION
### 4.1 Comparison between airfoil and P1 pipe diffuser at 100 % rotating speed

The performance map of optimized airfoil diffuser and baseline design of pipe diffuser P1 at 100% speed is as shown in Fig. 8. Isentropic efficiency and total-to-total pressure ratio of compressor stage were considered for comparison. At designed mass flow rate (0.806 m/s), airfoil diffuser and P1 pipe diffuser has isentropic efficiency as 84.05% and 81.62%, and total-to-total pressure ration as 4.42 and 4.33, respectively—representing 2.89% decrement in $\eta_{t-t}$ and 2.04\% decrement in $P_{t_4}/P_{t_1}$. Other than that, operating range of mass flow rate is also found to be wider in case of airfoil diffuser than with pipe diffuse; but the superiority in $\eta_{t-t}$ of airfoil over pipe reduces toward the surge limit, Fig. 8. Fig. 9 shows variation of $P_s$ and $P_t$ along streamline in pipe and airfoil diffuser. Pressure rise trend in both pipe and airfoil diffuser agrees well with the theroy of subsonic diffusion of compressible flow. From isentropic compressible flow equation, one could conlclude that the increase in dp/dx is more when $Ma$ is high [14]. Hence, it is expected that more pressure rise would occur in the front part. The red triangle marker indicates the location of throat. Since throat length is kept as minimum as possible (almost avoided), no local flow acceleration and static pressure drop after throat was found. Almost 60% of total pressure rise was observed before throat in case with pipe diffuser; However, the magnitude of total pressure rise is still lower than airfoil diffuser. It is interesting to note that despite the lower streamwise length from throat to outlet in pipe diffuser, almost 6% of $P_{t_{th}}$ drop was observed, which is 4% in case with airfoil diffuser. The throat blockage in pipe diffuser was also compared with airfoil diffuser using mass-averaged quantities at design point, refer **Appendix B**. The value of $B_{th}$ for airfoil and pipe diffuser is 0.16 and 0.02, respectively. Which indicates that blockage due to boundary layer thickness is much lesser in pipe diffuser than in airfoil diffuser—provides an opportunity to reduce the size of throat further. This can also be seen in Fig. 10. To gain further insight on poor performance on P1 pipe diffuser, discussion on detailed flow mechanism was made in subsequent paragraphs.

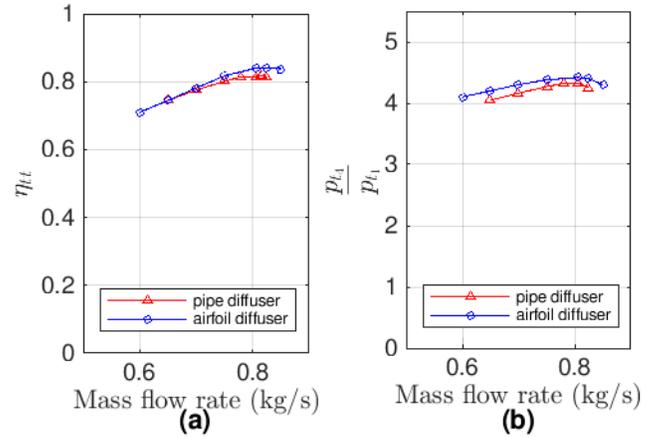

**FIGURE 8:** (a) Isentropic efficiency and (b) total-total pressure ratio

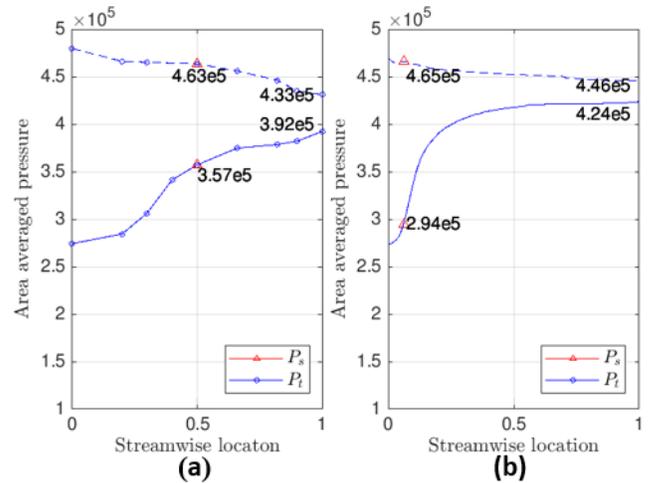

**FIGURE 9:** $P_s$ and $P_t$ variation with streamwise location in (a) pipe and (b) airfoil diffuser at design mass flow rate

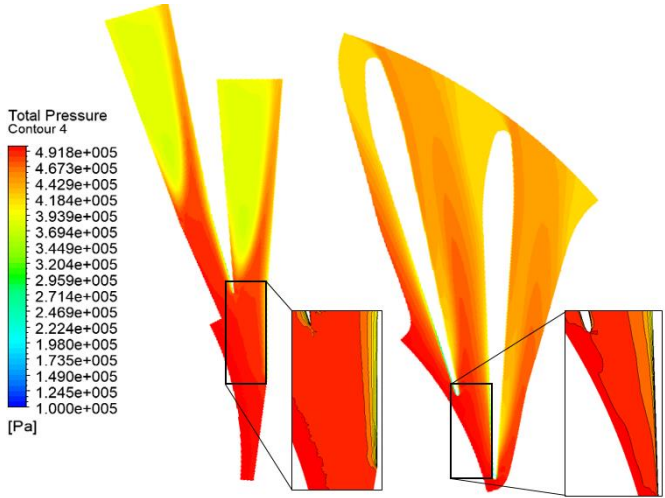

**FIGURE 10:** Total pressure contour at 50% span at design mass flow rate

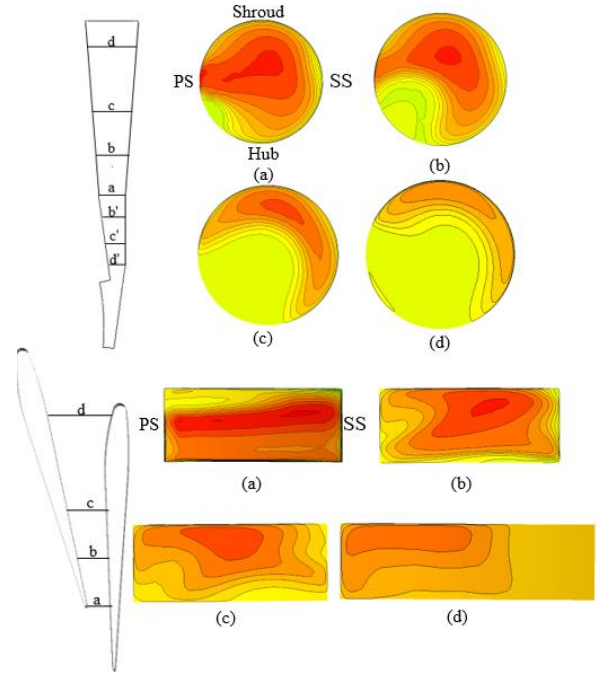

**FIGURE 11:** Total pressure contours at different locations along streamline

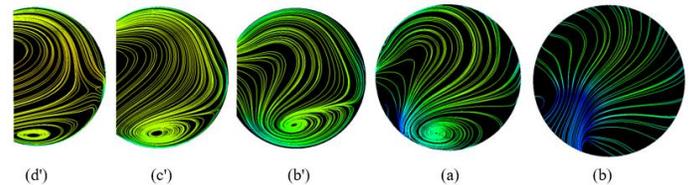

**FIGURE 12:** Vorticity development at leading edge of pipe diffuser

In Fig. 10, by observing the contours from leading edge to throat, one may find that the growth of boundary layer is more in airfoil diffuser, and more loss due to boundary layer thickness occurs in airfoil than in pipe diffuser. Total pressure contour at location (a) in Fig 11 also justifies the same. This observation is well supported by the findings of Kenny [1], Dean (1973) [6], and Ge Han et al. [5]. The ridges of the leading edge of pipe diffuser generate two counter-rotating vortices on the PS, also refer streamline plot at location (d') in Fig. 12. This would bring the low momentum fluid at near wall to a high momentum fluid at center, which leads to a better flow mixture in the pseudo- and semi-vaneless space. This phenomenon causes lower value of $B_{th}$, as discussed above.

The next question would be whether this counter-rotating vorticies carry the benefit further downstream into Channel space. It is found that under given operating parameters of MGT and for P1 pipe diffuser, it actually destabilizes the flow further downstream. Fig. 12 shows the development and mixing out of two-counter rotating vortices. The strength of two vortices is different, which might be because of non-uniformity of axial velocity leaving the impeller. It is found that vortices couldn't sustain and mixes out right after the throat, Fig 12. By considering the PS boundary layer, SS boundary layer, and two counter-rotating vortices; it is found that latter mechanism dominates the location of PS or SS separation. Zachau [7] found from experimental data that high level of vortex strength could initiate PS separation. This PS separation has different mechanism of separation than the separation due to boundary layer growth. The two counter-rotating vortices have tendency to bring low momentum fluid to PS and high-momentum fluid to SS, greater the strength of vortices larger the magnitude of this phenomena. It is evident by looking at Fig. 10 to Fig. 13 that something similar is happening here, which causes a large amount of PS separation further downstream in case with pipe diffuser. By looking at total pressure contours in airfoil diffuser in Fig. 10 and Fig. 11, it found that loss initiates at walls and become prominent at SS further downstream. Hence, from above discussion it is concluded that pipe diffuser exhibits PS separation while airfoil diffuser exhibits SS separation; However, the magnitude of former seems to be quite higher than later.

### 4.2 Span-wise characteristics of vorticity induced flow separation at design mass flow rate

Vorticity induced PS separation was discussed in previous section. To investigate the flow separation further downstream into the Channel space, five planes parallel to Shroud and Hub was created at different spans, as shown in Fig. 14. To simultaneously understand the state of streamline in orthogonal plane, four locations (a, b, c, d) were marked as shown in Fig. 14. From looking at the streamline plot at orthogonal planes at locations a, b, c and d; one could conclude that the development, or rather mixing out, of vortices is not symmetrical about 50% span plane. Flow separation is more likely to happen near Hub side than Shroud for all four locations. 10% span and 30% span fall in to the domain of low momentum wake flow cells for all four locations, so the flow is easier to separate as shown in Fig. 14. For 50% span, at location a and b, no separation occurs; it is evident from the streamline view in orthogonal plane that vortices rather help to reattach the flow near the PS wall at location a and b. However, further downstream at locations c and d, flow separations are likely to occur at PS. For 70% span, vertices benefit at locations a, b, and c, and causes PS separation at location d. For 90% span, it seems that vertices benefit in all four location. In summary, it is obvious that vortices have dual effect on flow separation. According to Fig. 14, the vortices induces flow separation near the vicinity of Hub, like 10%, 30%, 50% span. Meanwhile, it helps to reattach the flow near the vicinity of Shroud, like 50%, 70%, 90% span. More detail note on vorticity generation and its effect on flow separation can be found in Zachau et al [7] and Z.Sun et al [15].

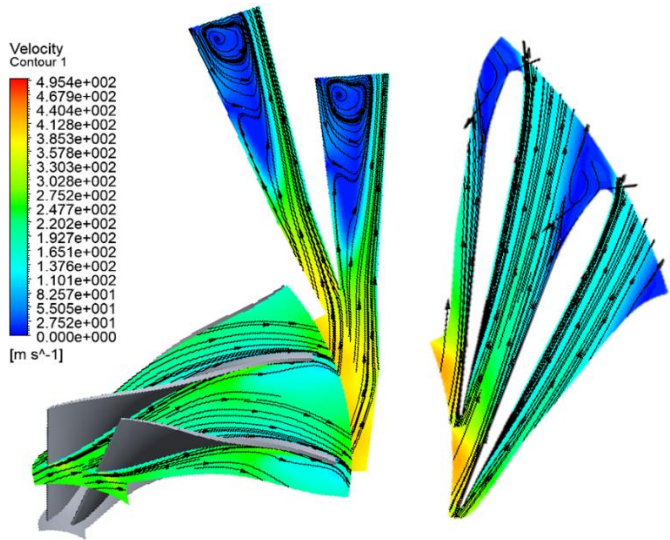

**FIGURE 13:** Velocity contour and vector at 50% span at design mass flow rate

### 4.3 Performance study of P1 pipe diffuser at selected off-design point

Performance map was generated at different mass flow rates for 100%, 90%, 80%, and 70% rotating speed, as shown in Fig. 15 and in Fig. 16, which shows plots of isentropic efficiency and total to total pressure ratio, respectively. In order to understand the performance at off-design points, detailed analysis was presented for four operating point—indicated as OP1, OP2, OP3, and OP4, Fig. 15. The first point, OP1, is the compressor design point with mass flow rate equals to 0.806 Kg/s. OP2 corresponds to near surge point at 100% rotating speed with mass flow rate equals to 0.65 Kg/s. Points OP3 and OP4 are design point and near surge point, respectively, for 60% rotating speed; with mass flow rate equals to 0.4 and 0.2 Kg/s, respectively.

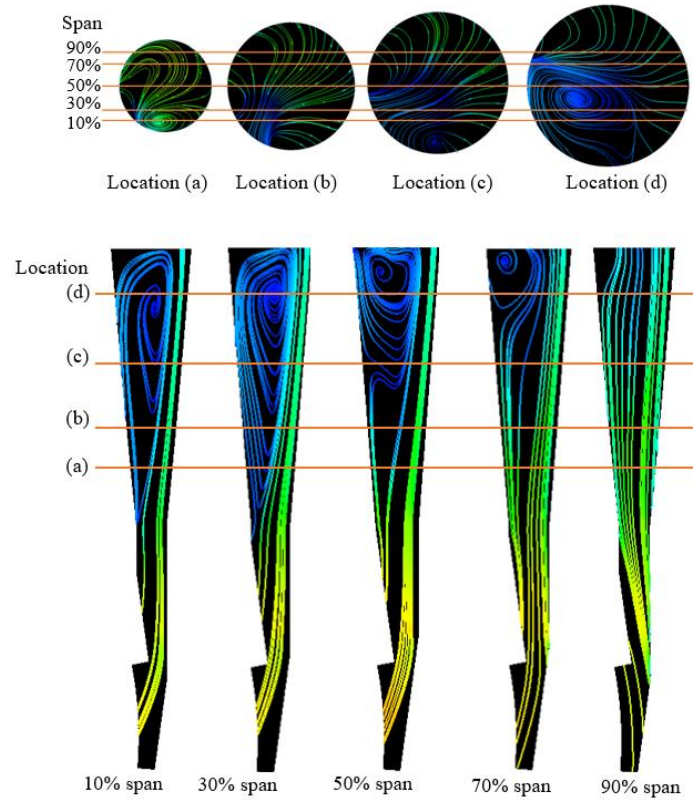

**FIGURE 14:** Flow separation at different spans

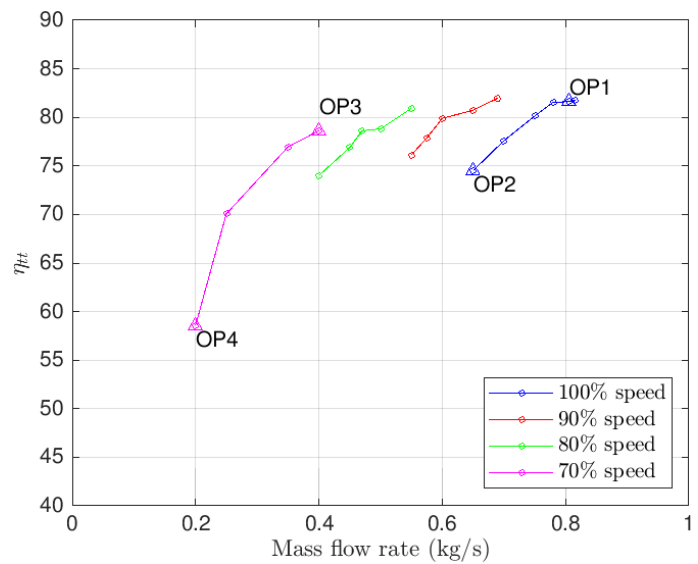

**FIGURE 15:** Isentropic efficiency plot at different operating speed

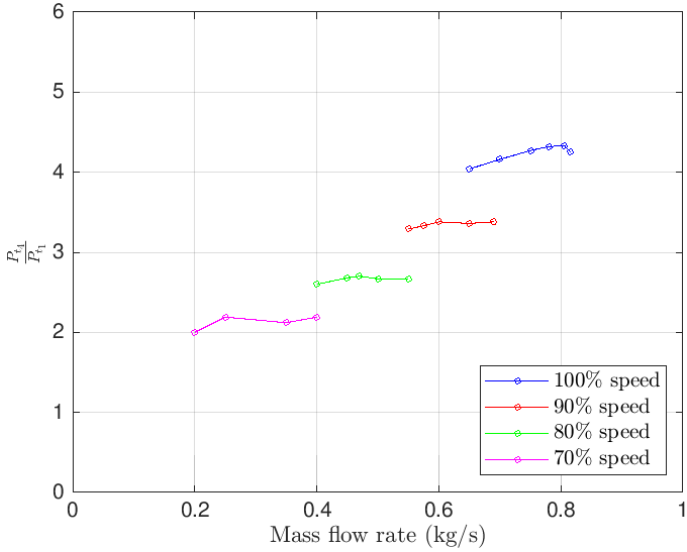

**FIGURE 16:** Total-to-total pressure ratio plot at different operating points

Figure 17 shows vorticiy at three orthogonal planes— locations (d'), (c'), and (a), for four operating points. Velocity contours and total pressure contours at 50% span plane for four operating points is represented as shown in Fig. 18. The Range of values of velocity and total pressure for contour generation was manually specified and was kept same across all operating points. Fig. 19 shows the pressure recovery and pressure loss coefficient variation along the streamline.

By observing Fig. 17 for OP2. it seems that behavior of vorticity development at leading edge is same as in case of OP1. However, at throat, no similarity exists. For OP2, due to low mass flow rate, inlet flow angle to diffuser increases. The increased inlet flow angle could increase throat blockage which might choke the flow at mass flow rate lower than expected. By observing total pressure contour (right) for OP2 in Fig. 18, it is evident that SS boundary layer have grown extensively due to misalignment of flow and geometrical angle. It can be seen that approximately 50\% of throat area is blocked by boundary layer, which leads to increased throat blockage. In case of OP2, growth of SS boundary layer dominates the mechanism of flow separation leading to SS separation, as shown in velocity contour (left) for OP2 in Fig. 18. Vorticity contour for OP2 at throat— location (a), in Fig. 17, indicates the low momentum wake flow cell zone near SS, which essentially means the same that out of vorticity induced PS separation and Boundary layer induced SS separation, later dominates in this case. SS separation has drastic effect on $C_p$ reduction than PS separation as can be inferred from the results of Wilkosz et al. [8]. The $C_p$ value in OP2 case comes out to be 0.47 as shown in Fig. 19. The value of loss coefficient increases to 0.32. By observing the trend of $C_p$ and $C_{p_o}$ it is clear that vanless and semi-vanless space offers no $C_p$ rise, rather increases the value of $C_{p_o}$. The trend is as expected.

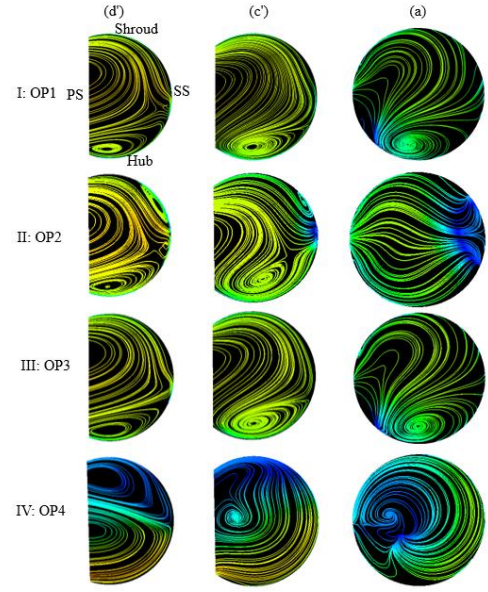

**FIGURE 17:** Isentropic efficiency plot at different operating speed

For OP3, the behavior of vorticity, velocity contour, and total pressure contour are essentially the same as OP1, Fig. 17 and Fig. 18. The trend of $C_p$ and $C_{p_o}$ variation along streamline is also quite matching with OP1 trend. The value of $C_p$ decreases very slightly to 0.55, and vale of $C_{p_o}$ increases slightly to 0.24, Fig. 19.

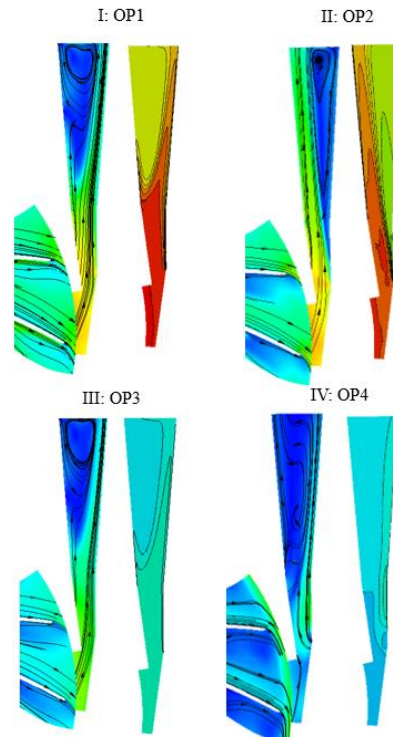

**FIGURE 18:** Velocity contour (left) and total pressure contour (right) at 50% span plane for four operating points

For OP4, by observing at Fig. 17, it seems that low momentum flow zone cluster at Shroud right from the inlet portion. However, the behaviour at throat is not inferential to author. By looking at total pressure contour in Fig. 18, it seems that growth of SS boundary layer is significant. But velocity contour shows PS separation further downstream. From observing $C_{po}$ trend in Fig. 19, it seems that more % of total pressure loss occur at inler portion than in Channel space. The values of $C_{po}$ increases to 0.55, and value of $C_p$ decreases to 0.42, Fig. 19.

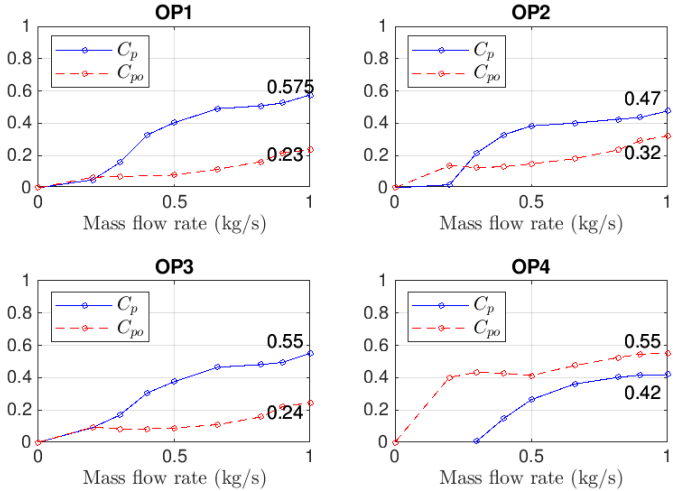

**FIGURE 19:** $C_p$ and $C_{po}$ variation along streamline

### 4.4 Parametric study on pipe diffuser

It is apparent that any effort to reduce the strength of vortices or to mitigate the vorticity induced PS separation would increase the performance of pipe over airfoil diffuser. In order to check the effect of various parameters, authors have tried to employ the statistical method of DOE module in ANSYS, but it has been observed that the atomization of meshing step was not accurately done, which lead to failure of this methodology. Hence, focus was made at conventional method of experimental design—OFAT. And few possible variations have been checked.

**TABLE 3:** Performance at different parameters at design point

| Designs | P1 | P2 | P3 | P4 |
|---|---|---|---|---|
| $\eta_{t-t}$ | 81.62% | 82.22% | 80.99% | 81.9% |
| $P_{s_4}$ | 3.93 | 4.10 | 3.62 | 3.93 |
| $P_{t_4}$ | 4.33 | 4.36 | 4.29 | 4.33 |

By decreasing $R_{tan}$, span of psuedo- and semi-vanless space can be reduced. Since, it works as vortex generator in pipe diffuser, it was predicted that a little decrease in its span might reduce the growth of vortices before throat and help to sustain it a little longer further downstream of throat. Though, there is not much difference in performance, a marginal improvement is found as shown in Table 3. It is also found that P2 design has $B_{th}$ equal to 0.06—slightly higher than P1, which is also as expected. Fig. 20 shows that slightly more total pressure loss and less pressure recovery occurs before throat in P2 as compared to P1; But the trend reverses after the throat and higher-pressure recovery of 0.607 could be achieved in P2 design. Purpose of P3 design is just to check whether the assumption of making $L$ similar to airfoil is convincing or not. In addition to that, results also support the inference from diffuser map [10] that 6° divergence is fine for given condition. Ge Han et al [2] avoided the PS separation by reducing $R_{3a}/R_3$ from 1.05 to 1.03. It was found that reducing $R_{3a}/R_3$ offers little advantage over P1, Table 3. However, none of the parameters have shown significant reduction in PS separation.

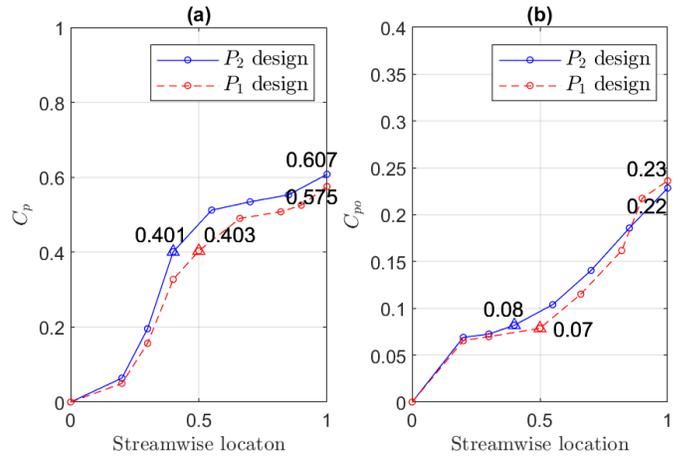

**FIGURE 20:** $C_p$ and $C_{po}$ variation along streamline for P1 and P2 designs

## 5 CONCLUSIONS
1. Under given dimensions of compressor impeller, pipe diffuser geometry offers little flexibility in independently deciding certain parameters such as $D_{th}$ and $\alpha$.
2. Leading ridges of pipe diffuser generates two counter rotating vortices of unequal strength. Which offers advantage in pseudo- and semi-vanless space. But destabilize the flow in Channel space leading to PS separation.
3. Throat blockage in pipe diffuser is found out to be less than airfoil diffuser. This is attributed to the presence of two counter-rotating vortices.
4. Interaction and transport of two vortices in Channel space help to reattach the flow at PS near Shroud region. But the separation at middle and near Hub region outweighs this effect, leading to inferior overall performance
5. Since higher the strength of vortices larger would be the PS separation, any attempt to reduce the strength under give operating parameter of micro gas turbine could increase its performance over airfoil diffuser. Improvement in geometry of leading ridge could be sought to balance the two-opposite effect of counter-rotating vortices.

6. $D_{tan}$ equals to 145 mm (which also equals to impeller outlet diameter $D_2$) has better performance than larger value of $D_{tan}$ such as 150 mm.

**Appendix A: Correlation of constrained parameters**

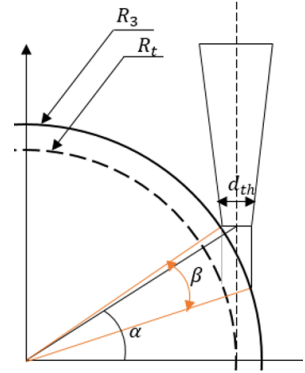

**FIGURE 21:** Scheme of pipe diffuser

By taking an assumption that throat length is zero (minimum value of throat length), circle with radius $R_3$ can be drawn as shown in Fig. $R_t$ is the radius of tangential reference circle, Fig. 20. This would give the maximum value of $N$. Exceeding this value would result in intersection of pipes beyond throat and geometry would fail. By performing little mathematical exercise on above schematic, authors have obtained following correlation ($R_t$ in this case refers to $R_{tan}$):

$$\beta = tan^{-1}\left(\frac{R_t\, tan\alpha}{R_t - 0.5D_t}\right) - cos^{-1}\left(\frac{R_t + 0.5D_t}{\sqrt{(R_t - 0.5D_t)^2 + (R_t\, tan\alpha)^2}}\right)$$

$$N_{max} = \left\lfloor \frac{\beta}{2\pi} \right\rfloor$$

**Appendix B: Throat blockage**

From isentropic compressible flow equation, area of isentropic flow or effective area ($A_{eff}$) can be calculated as:

$$A_{eff} = \left(\frac{m\sqrt{T_0}}{P_0}\right)\left(\frac{R}{\gamma}\right)\left(\frac{\left(1 + \frac{\gamma-1}{2}M^2\right)^{\frac{\gamma+1}{2(\gamma-1)}}}{M}\right)$$

$$1 - B_{th} = \frac{A_{eff}}{A_{geo}}$$